\documentclass[11pt,leqeq]{article}
\usepackage{amssymb,amsthm}
\usepackage{amsmath}
\usepackage{amscd}
\def\C{{\mathbb C}}

\def\ob{{\rm{Ob}}}

\def\C{{\cal{C}}}

\def\D{{\cal{D}}}

\def\FGD{{{gpd}}}

\def\Iso{{\rm{Iso}}}

\newtheorem{thm}{Theorem}%[section] (If you want theorem numbered
\newtheorem{cor}[thm]{Corollary}%       goes for lemmas, etc.)
\newtheorem{prop}[thm]{Proposition} %--> \begin\end{theorem,lemma,...}

\newtheorem{defi}[thm]{Definition}
\newtheorem{exa}[thm]{Example}

\newcommand{\des}{\displaystyle}
\begin{document}

\title{Rational Combinatorics}
\author{H\'{e}ctor Bland\'{i}n\
and Rafael D\'{i}az } \maketitle
\begin{abstract}
\noindent We propose a categorical setting for the study of the
combinatorics of rational numbers. We find combinatorial
interpretation for Bernoulli and Euler numbers and polynomials.\\

\end{abstract}
\section{Introduction}
\noindent Let $Cat$ be the category whose objects are small
categories (categories whose collection of objects are sets) and
whose morphisms $Cat(\C,\D)$ from category $\C$ to category $\D$
are functors $F:\C\rightarrow\D.$ Let $Set$  be the category of
sets. We define an equivalence relation $Iso_{\C}$ on $Ob(\C)$,
objects of $\C$, as follows:  $x$ and $y$ are equivalent if and
only if there exists an isomorphism $f\in\C(x,y).$ There is a
natural functor $D:Cat\rightarrow Set$ called
\textit{decategorification} given by
$D(\C)=Ob(\C)/{\Iso}_{\C}$ for any small category $\C$. Given a
functor $F:\C\rightarrow\D$ then
$D(F):Ob(\C)/{Iso_{\C}}\rightarrow Ob(\D)/{Iso_{\D}}$ is the
induced map. If $x$ is a set and $\C$ is a category such that
$D(\C)=x,$ we say that $\C$ is a
\textit{categorification} of $x.$ The reader may consult
\cite{BaezDolan3} and \cite{RDEP} for more on the notion of
categorification. \\

\noindent We are interested in the categorification of
sets with additional properties, for example one would like to
find out what is the categorification of a ring. It turns out that
the definition of categorification  given above is too rigid, for
most applications a weaker notion  seems to be more useful. For
example, see Section
\ref{catofrings} for details, a  categorification of a
ring $R$ is a category $\C$ together with bifunctors $\oplus$ and
$\otimes$ from $\C \times \C$ to $\C$, distinguished objects $0$
and $1$, and a negative functor $N$ from $\C$ to $\C$. Moreover
$\C$ should be provided with a valuation  map $|\ \
|:Ob(\C)\rightarrow R$ such that $|a|=|b|$ if $a$ and $b$ are
isomorphic, $|a
\oplus b|=|a|
\oplus |b|,$\ $|a\otimes b|=|a||b|,$\ $|1|=1$ and $|0|=0$, and
$|N(a)|=-|a|$, for any objects $a$ and $b$ in $\C$.
Categorification of
semi-rings is defined similarly but the functor $N$ is no longer required.\\

\noindent Our starting point is the theory of combinatorial species
of Joyal, see \cite{j1} and \cite{j2}, which can be described as
starting from a categorification of the natural numbers
$\mathbb{N}$ and extending it to a categorification of the
semi-ring $\mathbb{N}[[x]]$ of formal power series with
coefficients in $\mathbb{N.}$ Namely, the category $set$ of finite
sets with disjoint union and Cartesian product, together with the
map $|\
\ |:Ob(set)\rightarrow
\mathbb{N}$ sending a finite set $x$ to its cardinality $|x|$,
is a categorification of the natural numbers. Joyal goes on and
shows that the category $set^{\mathbb{B}}$ of functors from
$\mathbb{B}$ to $set$ defines a categorification of
the semi-ring $\mathbb{N}[[x]]$.\\

\noindent As explained by Zeilberger in \cite{DZeilberger} the main topic of
enumerative natural combinatorics is the following: given a
infinite sequence $x_{0},x_{1},\dots,x_{n},\dots$ of finite sets,
objects of $set,$ compute the corresponding sequence
$|x_{0}|,|x_{1}|,\dots,|x_{n}|,\dots$ of natural numbers. So one
looks for a numerical representation of combinatorial objects.
There is also an inverse problem in enumerative combinatorics:
given a sequence of natural numbers
$a_{0},a_{1},\dots,a_{n},\dots$ find an appropriated sequence of
finite sets $x_{0},x_{1},\dots,x_{n},\dots$ such that
$|x_{i}|=a_{i}$ for $i\in\mathbb{N}.$ Here we look for a
combinatorial interpretation of a sequence of natural numbers. In
the theory of species a fundamental part is played by the groupoid
$\mathbb{B}$ whose objects are finite sets and whose morphisms are
bijections. From the point of view of the theory of species the
main problem of enumerative natural combinatorics can be described
as follows: given a functor $F:\mathbb{B}
\rightarrow set$ find the sequence
$|F([0])|,|F([1])|,....,|F([n])|,...$ or what is essentially the
same, find the associated generating series
$\sum_{n=0}^{\infty}|F([n])|\frac{x^n}{n!}$.\\

\noindent Enumerative combinatorics can be extended to deal with integer
numbers. The main problem of enumerative integral combinatorics is
the following: given a sequence
$(x_{0},y_{0}),\dots,(x_{n},y_{n}),\dots$ of pairs of finite sets,
objects of  $\mathbb{Z}_{2}$-$set$ to be defined in Section
\ref{NRS},  compute the
corresponding sequence
$|(x_{0},y_{0})|,\dots,|(x_{n},y_{n})|,\dots$ of integers, where
for each pair of finite sets $(x,y)$ its cardinality is defined by
$|(x,y)|=|x| - |y|$. The inverse problem is the following: given a
sequence $a_0,\dots,a_n,\dots$ of integers find a nice sequence
$(x_{0},y_{0}),\dots,(x_{n},y_{n}),\dots$ of pairs of finite sets
such that $|(x_{i},y_{i})|:=|x_{i}|-|y_{i}|=a_{i}$ for
$i\in\mathbb{N}.$ From the species point of view the main problem
of integral combinatorics may be described as follows: given a
functor $F: \mathbb{B} \rightarrow  \mathbb{Z}_{2}$-$set$ compute
the associated generating series
$\sum_{n=0}^{\infty}|F([n])|\frac{x^n}{n!}$.\\

\noindent  We face the following problem in this paper: find a categorification
of $\mathbb{Q}$ the ring of rational numbers. In the first two
sections of this paper we concentrate on the problem of the
categorification of $\mathbb{Q}_{\geq 0}$, the semi-ring of
nonnegative rational numbers, leaving the study of the
categorification of $\mathbb{Q}$ for the final two sections. So we
need a category $\C$ provided with sum and product bifunctors
together with a valuation map $|\
\ |:Ob(\C)\rightarrow\mathbb{Q}_{\geq 0 }$ satisfying a natural
set of axioms.
%We argue that the category $\FGD$ of finite groupoids plays the
%part of $\mathbb{B}$ in the study of the combinatorial properties
%of rational numbers.
Given such pair $\C$ and $|\ \ |$ one defines the main problem of
"rational combinatorics" as the problem of finding the sequence
$|x_{0}|,\dots,|x_{n}|,\dots$ for any sequence
$x_{0},\dots,x_{n},\dots$ of objects of $\C.$ Similarly the
inverse problem in "rational combinatorics" would be the
following: given a sequence $a_{0},\dots,a_{n}\dots$ in
$\mathbb{Q}_{\geq 0},$ find a nice sequence
$x_{0},\dots,x_{n},\dots$ of objects of $\C$ such that
$|x_{i}|=a_{i},$  for $i\in\mathbb{N}.$ There are several
categories $\C$ provided with a valuation map $|\ \
|:Ob(\C)\rightarrow\mathbb{Q}_{\geq 0}$, so if we like our
enumerative problem above to be consider as being "combinatorial"
we should demand that the category $\C$  be close to $set$, and
the valuation map $|\ \ |$ close to the notion of
cardinality of finite sets. \\

\noindent Fortunately, Baez and Dolan in \cite{BaezDolan} have proposed a good
candidate to play the part of $set$ when dealing with the
combinatorial properties of rational numbers, namely the category
$\FGD$ of finite groupoids. A groupoid is a category such that all
its morphisms are invertible. According to Baez and Dolan the
cardinality of a finite groupoid ${G}$ is given by the formula
$|{G}|=
\sum_{x
\in D({G})}\frac{1}{|{G}(x,x)|}$. Starting
from this definition we construct a categorification of the ring
$\mathbb{Q}$ of rational numbers. This construction can be further
generalized using Joyal's theory of species to yield a
categorification of the ring $\mathbb{Q}[[x_1,...,x_n]]$ of formal
power series with rational coefficients in $n$-variables . This
paper is devoted to highlight the properties  that make $\FGD$
into a good ground set for the study of the combinatorial
properties of rational numbers, however we mention from the outset
its main shortcoming: it does not seem to exist a functorial way
to associate to each groupoid $G$ another groupoid $G^{-1}$ such
that $|G^{-1}|=|G|^{-1}.$

\section{Categorification of rings}\label{catofrings}
\noindent In this section we introduce the notion of categorification of rings.
A \textit{categorification} of a ring $R$ is a triple $(\C,N,|\ \
|)$ such that $\C$ is a category, $N:\C \rightarrow \C$ is a
functor and $|\ \ |: Ob(\C)
\rightarrow R$ is a map called the valuation map. This data is
required to satisfy the following axioms:

\begin{enumerate}
\item{$\C$ is provided with bifunctors $\oplus:\C \times \C
\rightarrow \C$ and $\otimes:\C \times \C
\rightarrow \C$ called sum and product, respectively.
Functors $\oplus$ and $\otimes$ are such that
\begin{itemize}
\item{There are distinguished objects $0$ and $1$ in $C$.}
\item{The triple $(\C,\oplus,0)$ is a symmetric monoidal category with
 unit $0$.}
\item{The triple $(\C, \otimes,1)$ is a  monoidal category
with unit $1$.}
\item{Distributivity  holds. That is  for objects $a,b,c$ of $\C$ there
are natural isomorphisms $a
\otimes(b\oplus c)\simeq (a \otimes b) \oplus  (a \otimes c)$ and $(a\oplus b)\otimes c \simeq (a \otimes c) \oplus (b \otimes c)$.}
\end{itemize}}
\item{The functor $N: \C \rightarrow \C$ must be such that for
objects $a,b$ of $\C$ the following identities hold
\begin{itemize}
\item{$N(a\oplus b)=N(a)\oplus N(b)$.}
\item{$N(0)=0$.}
\item{$N^{2}=I$ (identity functor).}
\end{itemize}}
\item{The map $|\ \ |:Ob(\C) \rightarrow R$ is such that for
objects $a,b$ of $\C$ the following identities hold
\begin{itemize}
\item{$|a|=|b|$ if $a$ and $b$ are isomorphic.}
\item{$|a\oplus b|=|a|+|b|$.}
\item{$|a \otimes b| = |a||b|$.}
\item{$|1|=1$ and $|0|=0$.}
\item{$|N(a)|=-|a|$.}
\end{itemize}}
\end{enumerate}

\noindent Let us make a few remarks regarding the notion of categorification
of rings.
\begin{enumerate}
\item{If $R$ is a semi-ring, i.e., we do not assume the existence of additive
inverses in $R$, then a categorification of $R$ is defined as
above but we omit the existence of the functor $N.$}
\item{$N$ is called the negative functor. In practice we prefer to write
$-a$ instead of $N(a)$.}
\item{Notice that we are not requiring that $\oplus$ and $\otimes$ be the coproduct and product of $\C$,
although in several examples they do agree.}

\item{See Laplaza's works \cite{l2} and \cite{l1} for the
full definition, and coherence theorems, of a category with two
monoidal structures satisfying the distributive property.}
\item{We stress that we only demand that $|a \oplus N(a)|=0$. Demanding the stronger condition
$a \oplus N(a)$ isomorphic to $0$, would reduce drastically the
scope of our definition.}
\item{We call $|a|$ the valuation of $a$. A categorification is
surjective if the valuation map is surjective.}
\item{$R$ is  a categorification of itself if we consider $R$ as
the category whose object set is $R$, and identities as morphisms.
The valuation map is the identity map, and the negative $N(r)$ of
$r \in R$  is $-r$. Thus, there is not existence problem related
with the notion of categorification: all rings admit a
categorification. The philosophy behind the notion of
categorification is that we can obtain valuable information about
a ring $R$ by looking at its various categorifications}.

\end{enumerate}

\noindent We close this section giving a couple of examples of
categorifications. Consider the category $set$ whose objects are
finite sets and whose morphisms are maps. The following data
describes $set$ as a surjective categorification of $\mathbb{N}$.
\begin{itemize}
\item{The triple $(set,\sqcup,\emptyset)$, where $\sqcup:set \times set \rightarrow
set$ is disjoint union, is a symmetric monoidal category with unit
$\emptyset$.}
\item{The triple $(set,\times,\{1\}),$ where $\times:set \times set \rightarrow
set$ is Cartesian product, is a monoidal category with unit
\{1\}}.
\item{The valuation  map $|\ \ |:Ob(set) \rightarrow \mathbb{N}$ sends a finite set $x$ to its cardinality $|x|$}.
\end{itemize}

\noindent Let $set^n$ be the category $set \times ....
\times set$. Objects in $set^n$ are pairs $(x,f)$ where $x$
is a finite set and $f:x \rightarrow \{1,...,n\}$ is a map.
Morphisms in $set^n$ from $(x,f)$ to $(y,g)$ are maps $\alpha: x
\rightarrow y$ such that $g \circ \alpha = f.$ Alternatively, objects
of $set^n$ can be describe as $n$-tuples $(x_1,...,x_n)$ of finite
sets.  A morphism in $set^{n}$ from $(x_1,...,x_n)$ to
$(y_1,...,y_n)$ is given by a $n$-tuple of maps
$(\alpha_1,...,\alpha_n)$ such that $\alpha_i:x_i \rightarrow y_i$
for $1\leq i \leq n.$ The following data describes $set^n$ as a
surjective categorification of $\mathbb{N}^{n}$.

\begin{itemize}
\item{The triple $(set^{n},\sqcup,\emptyset)$ is a  symmetric
monoidal category with unit $\emptyset=
(\emptyset,...,\emptyset).$ Disjoint union $\sqcup:set^{n}
\times set^{n} \rightarrow set^{n}$ is given by
$$(x_1,...,x_n)\sqcup(y_1,...,y_n)=(x_1 \sqcup y_1,...,x_n
\sqcup y_n).$$}
\item{The triple $(set^{n},\times,\{1\})$ is a monoidal
category with unit $\{1\}=(\{1\},...,\{1\})$.  Cartesian product
$\times:set^{n}
\times set^{n}
\rightarrow set^{n}$ is given by
$$(x_1,...,x_n)\times (y_1,...,y_n)=(x_1 \times y_1,...,x_n
\times y_n).$$}
\item{The valuation map $|\ \ |:Ob(set^{n}) \rightarrow \mathbb{N}^{n}$ sends  $(x_1,...,x_n)$
into  $(|x_1|,...,|x_n|)$}.
\end{itemize}

\noindent In the construction above we have implicitly made use of the
fact that if $\C_1$ and $\C_2$ are categorifications of rings
$R_1$ and $R_2$, respectively, then $\C_1
\times \C_2$ is a categorification of $R_1
\times R_2$, with  sum, product, negative functor and valuation map defined
componentwise.

\section{Nonnegative rational species}
We begin this section  introducing groupoids and their
cardinality, then we provide a list of useful examples of
groupoids. We use the following notations $[n]=\{1,2,\dots,n\},$
$S_{n}=\{f:[n]\rightarrow[n]\mid f\
\mathrm{is}\
\mathrm{bijective}\}$ is the symmetric group on $n$ letters,
$\mathbb{Z}_{n}=\{0,\dots,n-1\}$ is the cyclic group of order $n$.

\begin{defi}
A groupoid $G$ is a category such that all its morphisms are
invertible. We denote by $Gpd$ the full subcategory of $Cat$ whose
objects are groupoids.
\end{defi}

\noindent Let us introduce a few examples of groupoids.

\begin{exa}
Any category $\C$ has an underlying groupoid, which has the same
objects as $\C$ and whose morphisms are isomorphisms in $\C.$
\end{exa}

\begin{exa}
$\mathbb{B}$ denotes the groupoid whose objects are finite sets
and whose morphisms are bijections between finite sets.
$\mathbb{B}$ is the underlying groupoid of $set$.
\end{exa}

\begin{exa}
$\mathbb{B}^{n}$ denotes the groupoid whose objects are pairs
$(x,f),$ where $x$ is a finite set and $f:x\rightarrow[n]$ is a
map. Morphisms in $\mathbb{B}^{n}$ from $(x,f)$ to $(y,g)$ are
bijections $\alpha:x\rightarrow y$ such that $g\circ\alpha=f.$
$\mathbb{B}^{n}$ is the underlying groupoid of $set^{n}.$
\end{exa}

%$\mathbb{B}$ is a symmetric monoidal category as explained in
%\cite{RDEP} with monoidal structure given by $x\sqcup y$ (disjoint
%union) and $x\times y$ (cartesian product), for all
%$x,y\in\ob(\mathbb{B}).$

\begin{exa}
A group ${G}$ may be regarded as the groupoid $\overline{G}$ with
one object $1$ and $\overline{G}(1,1)=G.$
\end{exa}
\begin{defi}
A groupoid $G$ is called finite if its set of objects is finite,
and $G(x,y)$ is a finite set for $x,y$ objects of $G$. We denote
by $\FGD$ the full subcategory of $Gpd$ whose objects are finite
groupoids.
\end{defi}
\noindent  Disjoint
union and Cartesian product are given on $\FGD$ as restrictions of
the corresponding functors on $Cat$, see \cite{SMacLane}. The
disjoint union $\C\sqcup\D$, of categories $\C$ and $\D$, is the
category with objects $Ob(\C)\sqcup Ob(\D)$, and morphisms from
$x$ to $y$ given by
$$\C\sqcup\D(x,y)=\left\{\begin{array}{cc}
\C(x,y) & \mathrm{if}\ x,y\in Ob(\C),\\
\D(x,y) & \mathrm{if}\ x,y\in Ob(\D),\\
\emptyset    & \mathrm{otherwise}.
\end{array}\right.$$
%\mathrm{if}\ (x\in\ob(\C)\ \mathrm{and}\ y\in\ob(\D))\
%\mathrm{or}\ (x\in\ob(\D)\ \mathrm{and}\ y\in\ob(\C))

\noindent The Cartesian product $\C\times\D$, of categories $\C$ and $\D$, is the
category  such that $$Ob(\C\times
\D)=Ob(\C)\times Ob(\D),$$ and for
$(x_{1},y_{1}),(x_{2},y_{2})\in Ob(\C\times\D)$ we have
$\C\times\D((x_{1},y_{1}),(x_{2},y_{2}))=\C(x_{1},x_{2})\times\D(y_{1},y_{2}).$
The unit for disjoint union is the category $\emptyset$ with no
objects. A unit for Cartesian product is the category $1$ with one
object and
one morphism.\\

\noindent Let $\mathbb{Q}_{ \geq 0}$ be the semi-ring of nonnegative
rational numbers. Recall that $D(G)$ denotes the set of
isomorphisms classes of objects in the groupoid $G$.

\begin{thm}\label{Qcategorification}
 The map $|\
\ |:Ob(\FGD)\rightarrow\mathbb{Q}_{ \geq 0}$ given by
$\displaystyle{|G|=\sum_{x\in D(G)}\frac{1}{|G(x,x)|}} $ is a
surjective $\mathbb{Q}_{\geq 0}$-valuation on $\FGD.$
\end{thm}
\begin{proof}
Let $G$ and $H$ be finite groupoids, then
\begin{eqnarray*}
|G\sqcup H|&=&\sum_{x\in D(G\sqcup H)}\frac{1}{|\left(G\sqcup
H\right)(x,x)|}
=\sum_{x\in D(G)}\frac{1}{|G(x,x)|}+\sum_{x\in D(H)}\frac{1}{|H(x,x)|}\\
&=&|G|+|H|,
\end{eqnarray*}
and
\begin{eqnarray*}
|G\times H|&=&\sum_{(x,y)\in D(G\times H)}\frac{1}{|G\times
H(x,y)|}=\sum_{(x,y)\in D(G)\times D(H)}\frac{1}{|G(x,x)||H(y,y)|}\\
&=&\left(\sum_{x\in D(G)}\frac{1}{|G(x,x)|}\right)
\left(\sum_{y\in D(H)}\frac{1}{H(y,y)}\right)=|G||H|.
\end{eqnarray*}
$|1|=1$ and $|\emptyset|=0$. For each
$\frac{a}{b}\in\mathbb{Q}_{\geq 0}$ the groupoid
${\mathbb{Z}_{b}}^{\sqcup a}$ satisfies
$\left|{\mathbb{Z}_{b}}^{\sqcup a}\right|=\frac{a}{b}.$
\end{proof}

\noindent Following Baez and Dolan we call $|G|$ the cardinality of the groupoid $G$.
The outcome of Theorem
\ref{Qcategorification} is that the pair
$(\FGD, |\ \ | )$ is a surjective categorification of the
semi-ring $\mathbb{Q}_{\geq 0}.$ Notice that any finite set $x$
may be regarded as the finite groupoid whose set of objects is $x$
and with identity morphisms only. We have an
inclusion $i:set \rightarrow \FGD$ such that $|x|=|i(x)|.$\\

\noindent Given categories $\C$ and $\D$ we let $\C^{\D}$ be the category
whose objects are functors $F:\D\rightarrow\C$.  Morphisms from
$F$ to $G$ in $\C^{\D}$ are natural transformations
$T:F\rightarrow G.$
\begin{defi} The category of
$\mathbb{B}$-$\FGD$ species is the category ${\FGD}^{\mathbb{B}}.$
An object $F:\mathbb{B}\rightarrow\FGD$ is called a
$\mathbb{B}$-$\FGD$ species or a (nonnegative) rational species.
The category of (nonnegative) rational species in $n$-variables is
${\FGD}^{\mathbb{B}^n}.$
\end{defi}

\noindent We denote by ${{\FGD}_{0}}^{\mathbb{B}^{n}}$ the full
subcategory of ${\FGD}^{\mathbb{B}^{n}}$whose objects are functors
$F:\mathbb{B}^{n}\rightarrow{\FGD}$ such that
$F(\emptyset)=\emptyset.$ Next we define operations on rational
species, which indeed apply as well to ${Gdp}^{\mathbb{B}^n}$.
These operations may be regarded as examples of a fairly general
construction given in
\cite{RDEP}. We let $\Pi[x]$ be the set of  partitions of
$x.$
\begin{defi}\label{operations}
Let $F,G\in Ob({\FGD}^{\mathbb{B}^{n}}),$ $G_{1},\dots,G_{n}\in
Ob({\FGD}_{0}^{\mathbb{B}^{n}})$ and $(x,f)$ be an object of
$\mathbb{B}^{n}.$ The following formulae define product,
composition and derivative of $\mathbb{B}^{n}$-$\FGD$ species
\begin{enumerate}
\item $(F+G)(x,f)=F(x,f)\sqcup G(x,f).$ \item
$(FG)(x,f)=\des{\bigsqcup_{\ \ \ a\sqcup b=x}F(a,f\mid_{a})\times
G(b,f\mid_{b})}.$ \item $(F\times G)(x,f)=F(x,f)\times G(x,f).$
%\item
%$F(G_{1},\dots,G_{n})(x,f)=\des{\bigsqcup_{\underset{g:\pi\rightarrow
%x}
%{\pi\in\Pi[x]}}F(\pi,g)\times\prod_{b\in\pi}G_{g(b)}(b,f\mid_{b})}.$
%\item $F(G_{1},\dots,G_{n})(x,f)=\des{\bigsqcup_{\underset{g:\pi\rightarrow
%[n]}
%{\pi\in\Pi[x]}}F(\pi,g)\times\prod_{b\in\pi}G_{g(b)}(b,f\mid_{b})}.$
%\item $F(G_{1},\dots,G_{n})(x,f)=\des{\bigsqcup_{\underset{\hat{f}:\pi\rightarrow
%[n]}
%{\pi\in\Pi[x]}}F(\pi,\hat{f})\times\prod_{b\in\pi}G_{\hat{f}(b)}(b,f\mid_{b})}.$
\item $F(G_{1},\dots,G_{n})(x,f)=\des{\bigsqcup_{\underset{g:\pi\rightarrow
[n]}
{\pi\in\Pi[x]}}F(\pi,g)\times\prod_{b\in\pi}G_{g(b)}(b,f\mid_{b})}.$
\item For $i\in[n]$,
$\partial_{i}:\FGD^{\mathbb{B}^{n}}\rightarrow\FGD^{\mathbb{B}^{n}}$
is the functor such that $\partial_{i}(F)$ is given by
$$\partial_{i}F(x,f)=F(x\sqcup\{*\},f\sqcup\{(*,i)\}).$$
\end{enumerate}
\end{defi}

\noindent For
$a=(a_{1},\dots,a_{n})\in\mathbb{N}^{n}$ we write
$[a]=([a_1],...,[a_n])$, $a!=a_{1}!\dots a_{n}!$ and
$x^{a}=x_{1}^{a_{1}}\dots x_{n}^{a_{n}}.$
\noindent The isomorphism class of an object $(x,f)$ in $\mathbb{B}^{n}$ is
given by the $n$-tuple
$(|f^{-1}({i})|)_{i=1}^{n}\in\mathbb{N}^{n}.$
\begin{thm}\label{Teor6}
The map $|\ \
|:Ob({\FGD}^{\mathbb{B}^{n}})\longrightarrow\mathbb{Q}_{\geq
0}[[x_{1},\dots,x_{n}]]$ given by
\begin{equation*}
|F|(x_{1},\dots,x_{n})=\sum_{a\in\mathbb{N}^{n}}|F([a])|\frac{x^{a}}{a!}
\end{equation*}
is a $\mathbb{Q}_{\geq 0}[[x_{1},\dots,x_{n}]]$-valuation on
${\FGD}^{\mathbb{B}^{n}}.$ Moreover $|F\times G|=|F|\times|G|$,
$|\partial_{i}F|=\partial_{i}|F|$ and\\
$$|F(G_{1},\dots,G_{n})|=|F|(|G_{1}|,\dots,|G_{n}|).$$
Above $\times$ is the Hadamard product of series and $G_i$ is
assumed to be in $\FGD_{0}^{\mathbb{B}^{n}}$ for $1\leq i \leq n.$
\end{thm}
\noindent The valuation $|F|$ of a rational species $F$ is called its generating
series. The outcome of Theorem \ref{Teor6} is that
$({\FGD}^{\mathbb{B}^{n}}, |\ \ |)$ is a surjective
categorification of the semi-ring $\mathbb{Q}_{\geq
0}[[x_{1},\dots,x_{n}]]$. We proceed to introduce a list of
examples of rational species computing in each case the
corresponding generating series.
\begin{exa}
For $1\leq i\leq n$ the singular species
$\mathrm{X}_{i}:\mathbb{B}^{n}\rightarrow\FGD$ is such that
\begin{equation*}
X_{i}(x,f)=\left\{\begin{array}{cc} 1 & if\ |x|=1\ \ and\ f(x)=i,\\
\emptyset & \ otherwise.\end{array}\right.
\end{equation*}
for $(x,f)$ in $\mathbb{B}^{n}.$ Clearly $|X_i|=x_i \in
\mathbb{Q}_{\geq 0}[[x_{1},\dots,x_{n}]].$
\end{exa}
\begin{exa}
Let $\mathrm{S}^{N}:\mathbb{B}\rightarrow\FGD$ be the species
sending a finite set $x$  to the groupoid $S^{N}(x)$ given by
$Ob(S^{N}(x))=\{x\}$
%$\left\{\begin{array}{cc}\{x\}& if x\neq\phi,\\ \phi &
%x=\phi.\end{array}\right.$
and $S^{N}(x)(x,x)=S_{|x|}^{N}.$ We have
$$\left|S^{N}\right|=\sum_{n=0}^{\infty}\frac{x^{n}}{(n!)^{N+1}}.$$
\end{exa}
\begin{exa}\label{Exa12}
Let $\mathbb{Z}^{N}:\mathbb{B}\rightarrow\FGD$ be such that for a
finite set $x$ the groupoid $\mathbb{Z}^{N}(x)$ is given by
$Ob(\mathbb{Z}^{N}(x))=\left\{\begin{array}{cc}
\{x\} & if x\neq\emptyset \\ \emptyset  & if x=\emptyset
\end{array}\right.$ and  $\mathbb{Z}^{N}(x)(x,x)=\mathbb{Z}_{|x|}^{N}$
for $x$ nonempty. We have
\begin{equation*}
\left|\mathbb{Z}^{N}\right|=\sum_{n=1}^{\infty}\frac{1}{n^{N}}\frac{x^{n}}{n!}.
%=\prod_{p\
%\mathrm{prime}}\left(1+\sum_{k=1}^{\infty}\frac{x^{p^{k}}}{(p^{k})!}p^{-kN}\right)
\end{equation*}
For $N=1$ we obtain
\begin{equation*}
|\mathbb{Z}|=
\sum_{n=1}^{\infty}\frac{1}{n}\frac{x^{n}}{n!}=\int\frac{e^{x}-1}{x}dx.
\end{equation*}
\end{exa}

\begin{exa}
Let $E_{N}:\mathbb{B}\rightarrow\FGD$ be such that for a finite
set $x$ the groupoid $E_{N}(x)$ is given by $Ob(E_{N}(x))=\{x\}$
and $E_{N}(x)(x,x)=\mathbb{Z}_{N}^{|x|}$. We have
\begin{equation*}
\left|E_{N}\right|=\sum_{n=0}^{\infty}\frac{1}{N^{n}}\frac{x^{n}}{n!}=e^{\frac{x}{N}}.
\end{equation*}
\end{exa}

\begin{exa}
\noindent Let $G$ be a group. Let
$\overline{G}:\mathbb{B}\rightarrow\FGD$ be the rational species
such that $\overline{G}(\emptyset)=\overline{G}$  and
$\overline{G}(x)=\emptyset$ if $x$ is nonempty. Clearly
$\left|\overline{G}\right|=\frac{1}{|G|}.$
\end{exa}

%\begin{exa}
%Let $\underline{\mathbb{Z}_{N}}:\mathbb{B}\rightarrow\FGD$ the
%rational species such that for $x\in\ob(\mathbb{B})$
%\begin{equation*}
%\ob\left(\underline{\mathbb{Z}_{N}}(x)\right)=
%\left\{
%\begin{array}{cc}
%\phi & if\ x\neq\phi, \\
%\mathbb{Z}_{N} & if\ x=\phi.
%\end{array}
%\right.
%\end{equation*}
%and Then we have
%\begin{equation*}
%\left|\underline{\mathbb{Z}_{N}}\right|=\frac{1}{N}.
%\end{equation*}
%\end{exa}
\begin{exa}
Let $(\mathbb{Z}_N )^{(\ \ )}:\mathbb{B}\rightarrow\FGD$ be the
rational species such that for a finite set $x$ the groupoid
$(\mathbb{Z}_N )^{(x)}$ is given by
\begin{equation*}
Ob\left((\mathbb{Z}_N
)^{(x)}\right)=\left\{\begin{array}{cc}\{x\}& if
x\neq\emptyset\\\emptyset & if\ x=\emptyset ,\end{array}\right.
\end{equation*}
and
\begin{equation*}
(\mathbb{Z}_{N}
)^{(x)}(x,x)=\mathbb{Z}_{N}\times\mathbb{Z}_{N+1}\times\dots\mathbb{Z}_{N+|x|-1}.
\end{equation*}
We have
\begin{equation*}
\left|(\mathbb{Z}_N  )^{(\ \ )}\right|=\sum_{n=1}^{\infty}\frac{1}{(N)^{(n)}}\frac{x^{n}}{n!}.
%=\frac{1}{N!}\int^{(N)}\left(\frac{e^{x}-\pi_{N}(e^{x})}{x^{N}/N!}\right)dx^{N}
\end{equation*}
\noindent where $(a)^{(b)}=a(a+1)\dots(a+b-1)$ is the increasing
factorial \cite{GCRota}.
\end{exa}
\begin{exa}
Let $P:\mathbb{B}\rightarrow\FGD$  be the rational species such
that $Ob(P(x))=\{a \mid a \subseteq x \}.$ Morphisms are given by
$P(x)(a,b)=\{\alpha:a\rightarrow b\mid\alpha\
\mathrm{is}\
\mathrm{bijective}\}$. We have
\begin{eqnarray*}
|P|=\sum_{n=0}^{\infty}\left(\sum_{k=0}^{n}\frac{1}{k!}\right)\frac{x^{n}}{n!}
%=\sum_{n=0}^{\infty}\left(\sum_{k=0}^{n}\frac{(n)_{k}}{k!}\frac{1}{(n)_{k}}\right)\frac{x^{n}}{n!}
%&=&\sum_{n=0}^{\infty}\left(\sum_{k=0}^{n}\frac{\left(\begin{array}{c}n\\
%k\end{array}\right)}{(n)_{k}}\right)\frac{x^{n}}{n!}.\\
%&=&\left(\sum_{n=0}^{\infty}\frac{1}{n!}\frac{x^{n}}{n!}\right)
%\left(\sum_{n=0}^{\infty}\frac{x^{n}}{n!}\right)
\end{eqnarray*}
\end{exa}
%\begin{defi}
%Let $f:D_{f}\subseteq\mathbb{R}\rightarrow\mathbb{R}$ and
%$g:D_{g}\subseteq\mathbb{R}\rightarrow\mathbb{R}$ a $N$ times
%differentiable function such that $\frac{d^{N}}{dx^{N}}g(x)=f(x)$
%we shall write $g(x)=\int^{(N)}f(x)dx.$
%\end{defi}
\noindent If $f$ and $g$ are formal power series such that
$\frac{d^{N}g(x)}{dx^{N}}=f(x),$ then we write
$\des{g(x)=\int^{(N)}}f(x)dx.$
\\

\noindent
%Let $\mathbb{Z}:\mathbb{B}\rightarrow\FGD$ the species
%given in example \ref{Exa12} with $N=1,$ and
%as follows for each $x\in\ob(\mathbb{B}),$ $\mathbb{Z}(x)$ is
%the groupoid such that
%$\ob(\mathbb{Z}(x))=\left\{\begin{array}{cc} \{x\} & if\ x\neq\phi,\\
%\phi & if\ x=\phi.\end{array}\right.$ and
%$\mathbb{Z}(x)(x,x)=\mathbb{Z}_{|x|}.$
\begin{defi}
The increasing factorial rational species
$\mathbb{Z}^{(N)}:\mathbb{B}\rightarrow\FGD$ is such that for each
finite set $x$ the groupoid $\mathbb{Z}^{(N)}(x)$ is given by
\begin{equation*}
Ob(\mathbb{Z}^{(N)}(x))=\left\{\begin{array}{cc}\{x\}&
x\neq\emptyset\\ \emptyset & if\ x=\emptyset, \end{array}\right.
\end{equation*}
and
$\mathbb{Z}^{(N)}(x)(x,x)=\mathbb{Z}_{|x|}\times\mathbb{Z}_{|x|+1}
\times\dots\times\mathbb{Z}_{|x|+N-1}$ for $x$ nonempty.
\end{defi}
\begin{thm}
\begin{equation*}
\left|\mathbb{Z}^{(N)}\right|=x^{1-N}\int^{(N)}{\frac{e^{x}-1}{x}}dx.
\end{equation*}
\end{thm}
\begin{proof}
%\begin{eqnarray*}
%|\mathbb{Z}^{(N)}|&=&\sum_{n=1}^{\infty}\frac{1}{(n)^{(N)}}\frac{x^{n}}{n!}
%=\frac{1}{x^{N-1}}\sum_{n=1}^{\infty}\frac{1}{n}\int^{(N-1)}x^{n}dx=
%\frac{1}{x^{N-1}}\int^{(N-1)}{\left(\sum_{n=1}^{\infty}\frac{x^{n}}{n}\right)}dx\\
%&=&\frac{1}{x^{N-1}}\int^{(N-1)}{\log\left(\frac{1}{1-x}\right)dx}
%={-x^{1-N}}\int^{(N-1)}\log(1-x)dx.
%\end{eqnarray*}
\begin{eqnarray*}
\left|\mathbb{Z}^{(N)}\right|&=&\sum_{n=1}^{\infty}\frac{1}{(n)^{(N)}}\frac{x^{n}}{n!}
=\frac{1}{x^{N-1}}\sum_{n=1}^{\infty}\frac{1}{n}\int^{(N-1)}\frac{x^{n}}{n!}dx=
\frac{1}{x^{N-1}}\int^{(N-1)}{\left(\sum_{n=1}^{\infty}\frac{1}{n}\frac{x^{n}}{n!}\right)}dx\\
&=&\frac{1}{x^{N-1}}\int^{(N-1)}{\left(\int{\frac{e^{x}-1}{x}}dx\right)dx}
={x^{1-N}}\int^{(N)}\frac{e^{x}-1}{x}dx.
\end{eqnarray*}
\end{proof}

\noindent Let $G$ be a group acting on a finite set $x.$ The
quotient groupoid $x/G$ is such that $Ob(x/G)=x$ and
$x/G(a,b)=\{g\in G\mid ga=b\}$ for  $a,b\in x.$ The
decategorification $D(x/G)$ of $x/G$ is just the quotient set of
the action of $G$ on $x$. Below we use the notation $O(a)= \{ga
\mid g
\in G\}.$
\begin{prop}
If $G$ acts on $x$ then $|x/G|=\frac{|x|}{|G|}.$
\end{prop}\begin{proof}\begin{equation*}
|x/G|=\sum_{\overline{a}\in D(x/G)}\frac{1}{|\rm{Iso}(\it
a)|}=\sum_{\overline{a}\in
D(x/G)}\frac{|O(a)|}{|G|}=\frac{1}{|G|}\sum_{\overline{a}\in
D(x/G)}|O(a)|=\frac{|x|}{|G|}.
\end{equation*}
\end{proof}
\noindent Let $G$ be a subgroup of $S_{k}.$
\begin{exa}
Let $P_{G}:\mathbb{B}\rightarrow\FGD$ be such that for a finite
set $x$ the groupoid  $P_{G}(x)$ is the quotient groupoid
$x^{k}/G.$ It should be clear that
$\des{|P_{S_{k}}|=\sum_{n=1}^{\infty}}\frac{n^{k}}{k!}\frac{x^{n}}{n!}$
and
$\des{|P_{\mathbb{Z}_{k}}}|=\sum_{n=1}^{\infty}\frac{n^{k}}{k}\frac{x^{n}}{n!}.$
\end{exa}

\noindent The inertia functor
$I:{\FGD}\rightarrow{\FGD}$ is such that for each $G$ the groupoid
$I(G)$ is given by
$$Ob(I(G))=Ob(G),$$ and for objects $a$ and $b$ of $I(G)$ we set
$$I(G)(a,b)=\left\{\begin{array}{cc} G(a,a) & \mathrm{if}\ \ a=b,
\\ \emptyset & \mathrm{otherwise}.\end{array} \right.$$

\noindent For example the set of objects of
$I\left(x^{k}/S_{k}\right)$ is $x^{k}.$ For $a,b\in x^{k}$ we have
\begin{equation*}
I\left(x^{k}/S_{k}\right)(a,b)=\left\{\begin{array}{cc}\{\sigma\in
S_{k}\mid a\sigma=a\} & \mathrm{if} \  a=b, \\ \phi & a\neq
b.\end{array}\right.
\end{equation*}
\noindent One checks that
\begin{equation*}
\left|\{\sigma\in S_{k}\mid a\sigma =a\}\right|= \prod_{i\in
x}|a^{-1}(i)|!.
\end{equation*}
\noindent Therefore if we assume that $|x|=n$ we get
\begin{equation*}
\left|I\left(x^{k}/S_{k}\right)\right|= \sum_{a:[k]\rightarrow
x}\frac{1}{\prod_{i\in X}|a^{-1}(i)|!} =\sum_{s_{1}+\dots+s_{n}=k}
\left(\begin{array}{c} k \\
s_{1},\dots,s_{n}\end{array}\right)\frac{1}{s_{1}!\dots s_{n}!}.
\end{equation*}

\noindent We extend the inertia functor to rational species
$I:{\FGD}^{\mathbb{B}} \rightarrow {\FGD}^{\mathbb{B}}$ by the
rule $I(F)(x)= I(F(x))$ for any finite set $x$. With this notation
we have shown the following result.
\begin{prop}
\begin{equation*}
\left|I\left(P_{S_{k}}\right)\right|=\sum_{n=0}^{\infty}\left(\sum_{s_{1}+\dots+s_{n}=k}
\frac{k}{(s_{1}!\dots s_{n}!)^{2}}\right)
\frac{x^{n}}{n!}.
\end{equation*}
\end{prop}

\noindent Let us introduce other interesting examples of rational species.
Let $\mathrm{Isinh}:\mathbb{B}\rightarrow\FGD$ be such that for a
finite set $x$ the groupoid $\mathrm{Isinh}(x)$ is given by
\begin{equation*}
Ob(\mathrm{Isinh}(x))=\left\{\begin{array}{cc}\{x\}& \mathrm{if}\
|x|\ \mathrm{is}\
\mathrm{odd},
\\ \emptyset\ &\ \mathrm{otherwise}.\end{array}\right.
\end{equation*}
\begin{equation*}
\mathrm{Isinh}(x)(x,x)=\mathbb{Z}_{|x|} \ \mathrm{for}\ |x|\  \mathrm{odd}.
\end{equation*}

\begin{prop}
\begin{equation*}
\des{|\mathrm{Isinh}|}=\des{\int\frac{\sinh(x)}{x}dx}.
\end{equation*}
\end{prop}
\begin{proof}
\begin{equation*}
\des{|\mathrm{Isinh}|=\sum_{n=0}^{\infty}\frac{1}{2n+1}\frac{x^{2n+1}}{(2n+1)!}}
=\des{\int\frac{\sinh(x)}{x}dx}.
\end{equation*}
\end{proof}
\noindent The rational species
$\mathrm{Icosh}:\mathbb{B}\rightarrow\FGD$ is such that for a
finite set $x$ the groupoid $\mathrm{Icosh}(x)$ is given by
\begin{equation*}
Ob(\mathrm{Icosh}(x))=\left\{\begin{array}{cc}\{x\}& \mathrm{if} \
|x|\ \mathrm{is}\
\mathrm{even},
\\ \emptyset\ &\ \mathrm{otherwise}.\end{array}\right.
\end{equation*}
\begin{equation*}
\mathrm{Icosh}(x)(x,x)=\mathbb{Z}_{|x|}\ \ \mathrm{for}\ x\ \mathrm{even}.
\end{equation*}

\begin{prop}
\begin{equation*}
\des{|\mathrm{Icosh}|}=\des{\int\frac{\cosh(x)}{x}dx}.
\end{equation*}
\end{prop}
\begin{proof}
\begin{equation*}
\des{|\mathrm{Icosh}|=\sum_{n=0}^{\infty}\frac{1}{2n}\frac{x^{2n}}{(2n)!}}=
\des{\int\frac{\cosh(x)}{x}dx}.
\end{equation*}
\end{proof}
\noindent Let us closed this section by  introducing a couple of species that will be used below.
\begin{exa}
The species $1\in{\FGD}^{\mathbb{B}^{n}}$ is such that for each
object $(x,f)$ of $\mathbb{B}^{n}$
$$1(x,f)=\left\{\begin{array}{cc}1 & \mathrm{if}\ x=\emptyset,\\
\emptyset &\ \mathrm{otherwise}. \end{array} \right.$$ Clearly $|1|=1.$
\end{exa}
\begin{exa}
The exponential species $\mathrm{Exp}\in\FGD^{\mathbb{B}^{n}}$ is
given by $\mathrm{Exp}(x,f)=1.$ Clearly we have
$\left|\mathrm{Exp}\right|=e^{x_{1}+\dots+x_{n}}.$
\end{exa}

\section{Rational species}\label{NRS}

\noindent So far we have consider only the categorification of
semi-rings. In order to proceed further and, in particular, make
sense of  multiplicative inverses in the category of species it
becomes necessary to deal with negative species. As explained in
Section \ref{catofrings} the categorification of a ring requires
in addition to the existence of sum, product, units and valuation,
the existence of a negative functor. In this section we are going
to construct categorifications of the rings $\mathbb{Z}$,
$\mathbb{Z}[[x_1,...,x_n]]$, $\mathbb{Q}$ and
$\mathbb{Q}[[x_1,...,x_n]]$, but we stress that our techniques
may be applied to many other rings as well.\\

\noindent We begin studying the categorification of $\mathbb{Z}$. Consider the
category $\mathbb{Z}_{2}\mbox{-}set=set \times set.$ Objects in
$\mathbb{Z}_{2}$-$set$ are pairs of finite sets. Morphisms from
$(a_1,b_1)$ to $(a_2,b_2)$, objects of $\mathbb{Z}_{2}$-$set$, are
given  by
$$\mathbb{Z}_{2}\mbox{-}set((a_1,b_1),(a_2,b_2))=set(a_1,a_2)
\times set(b_1,b_2).$$

\noindent The disjoint union bifunctor $\sqcup
:\mathbb{Z}_{2}\mbox{-}set
\times \mathbb{Z}_{2}\mbox{-}set \rightarrow \mathbb{Z}_{2}\mbox{-}set$ is
given by
$$(a_1,b_1)\sqcup(a_2,b_2)=(a_1 \sqcup a_2,b_1 \sqcup b_2).$$
The Cartesian product bifunctor is given by
$$(a_1,b_1)\times (a_2,b_2)=
(a_1 \times a_2 \ \sqcup \ b_1  \times b_2 \ , \ a_1 \times b_2 \
\sqcup
\  a_2
\times b_1).$$

\noindent The distinguished objects $\emptyset$ and $1$ are
$\emptyset = (\emptyset, \emptyset)$ and $1= (1,\emptyset)$,
respectively. The valuation map $|\ \ |:\mathbb{Z}_{2}\mbox{-}set
\rightarrow \mathbb{Z}$ is given by $|(a,b)|=|a|-|b|$, the negative functor
$N:\mathbb{Z}_{2}\mbox{-}set
\rightarrow \mathbb{Z}_{2}\mbox{-}set$  is given by $N(a,b)=(b,a)$
for  $(a,b)$ in $\mathbb{Z}_{2}\mbox{-}set.$ With these
definitions we obtain the following result.

\begin{thm}
$(\mathbb{Z}_{2}\mbox{-}set, N, |\
\ |)$  is a surjective categorification of $\mathbb{Z}.$
\end{thm}

\noindent The categorifications
$set$ of $\mathbb{N}$ and $\mathbb{Z}_{2}\mbox{-}set$ of
$\mathbb{Z}$ are fundamentally different. In $set$ objects $x$ and
$y$ are isomorphic if and only if $|x|=|y|.$ That is not the case
in $\mathbb{Z}_{2}\mbox{-}set,$ where for example
$|([2],[4])|=|([11],[13])|=-2$ but $([2],[4])$ and $([11],[13])$
are not isomorphic . At first sight this may seem like a nuance,
however, there is nothing wrong with the fact that $([2],[4])$ and
$([11],[13])$ are not isomorphic, since it is rather intuitive
that they provide different combinatorial interpretations for
$-2.$ Indeed the pair $([2],[4])$ leads to an interpretation of
$-2$ as a difference of even numbers, while the pair $([11],[13])$
leads to an interpretation of  $-2$ as a difference of prime
numbers.\\

\noindent We are ready to handle negative species.\\

\begin{defi} The category of
$\mathbb{B}$-$\mathbb{Z}_{2}\mbox{-}set$ species is the category
${\mathbb{Z}_{2}\mbox{-}set}^{\mathbb{B}}.$ An object
$F:\mathbb{B}\rightarrow \mathbb{Z}_{2}\mbox{-}set$ is called a
$\mathbb{B}$-$\mathbb{Z}_{2}\mbox{-}set$ species or an integral
species. The category of integral species in $n$-variables is
${\mathbb{Z}_{2}\mbox{-}set}^{\mathbb{B}^n}.$
\end{defi}

\noindent We define operations on integral species using the same formulae as
in Definition \ref{operations} but instead of using the monoidal
structures of $\FGD$ we use the corresponding structures in
$\mathbb{Z}_{2}\mbox{-}set.$ The negative functor
$N:\mathbb{Z}_{2}\mbox{-}set^{\mathbb{B}^n} \rightarrow
\mathbb{Z}_{2}\mbox{-}set^{\mathbb{B}^n}$ is defined as
$N(F)(x,f)=N(F(x,f))$ for $F$ in
$\mathbb{Z}_{2}\mbox{-}set^{\mathbb{B}^n}$. The valuation map $|\
\ |: Ob(\mathbb{Z}_{2}\mbox{-}set^{\mathbb{B}^n})
\rightarrow
\mathbb{Z}[[x_1,...,x_n]]$ is given by $|F|= \sum_{a \in
\mathbb{N}^{n}}|F([a])|\frac{x^{a}}{a!}$.

\begin{thm}
$(\mathbb{Z}_{2}\mbox{-}set^{\mathbb{B}^n}, N, |\ \ |)$ is a
surjective categorification of $\mathbb{Z}[[x_1,...,x_n]].$
\end{thm}

\noindent  We proceed to study the
categorification of the ring of rational numbers. Consider the
category $\mathbb{Z}_{2}\mbox{-}gpd=gpd \times gpd$ whose objects
are pairs of finite groupoids. Morphisms are given by
$$\mathbb{Z}_{2}\mbox{-}gpd((a_1,b_1),(a_2,b_2))=gpd(a_1,a_2)
\times gpd(b_1,b_2),$$
for  $(a_1,b_1),(a_2,b_2)$ objects of
$\mathbb{Z}_{2}\mbox{-}gpd$.\\

\noindent Disjoint union  $\sqcup :\mathbb{Z}_{2}\mbox{-}gpd
\times \mathbb{Z}_{2}\mbox{-}gpd \rightarrow \mathbb{Z}_{2}\mbox{-}gpd$ is
given by $(a_1,b_1)\sqcup(a_2,b_2)=(a_1 \sqcup a_2,b_1 \sqcup
b_2).$  Cartesian product is given by $(a_1,b_1)\times (a_2,b_2)=
(a_1 \times a_2 \ \sqcup \ b_1  \times b_2 \ , \ a_1 \times b_2 \
\sqcup
\  a_2
\times b_1).$ Distinguished objects  are  $\emptyset
= (\emptyset, \emptyset)$ and $1= (1,\emptyset)$. The valuation
map $|\ \ |:\mathbb{Z}_{2}\mbox{-}gpd
\rightarrow \mathbb{Z}$ is given by $|(a,b)|=|a|-|b|$, and the negative functor
$N:\mathbb{Z}_{2}\mbox{-}gpd
\rightarrow \mathbb{Z}_{2}\mbox{-}gpd$ is given by $N(a,b)=(b,a),$
for  $(a,b)$ in $\mathbb{Z}_{2}\mbox{-}gpd.$ We have the following
result.
\begin{thm}
$(\mathbb{Z}_{2}\mbox{-}gpd, N, |\
\ |)$  is a surjective categorification of
$\mathbb{Q}.$
\end{thm}

\noindent Objects of $\mathbb{Z}_{2}\mbox{-}gpd$ are called $\mathbb{Z}_{2}$-graded groupoids
or $\mathbb{Z}_{2}$-groupoids. We write $a-b$ instead of $(a,b)$
for  $(a,b)$ in $\mathbb{Z}_{2}\mbox{-}gpd$. For example, $-a$
denotes the $\mathbb{Z}_{2}$-groupoid $(\emptyset,a)$. We are
ready to give a full definition of the category of rational
species.

\begin{defi} The category of
$\mathbb{B}$-$\mathbb{Z}_{2}\mbox{-}gpd$ species is
${\mathbb{Z}_{2}\mbox{-}gpd}^{\mathbb{B}}.$ An object
$F:\mathbb{B}\rightarrow \mathbb{Z}_{2}\mbox{-}gpd$ is called a
$\mathbb{B}$-$\mathbb{Z}_{2}\mbox{-}gpd$ species or a rational
species. The category of rational species in $n$-variables is
${\mathbb{Z}_{2}\mbox{-}gpd}^{\mathbb{B}^n}.$
\end{defi}

\noindent Operations on rational species are constructed as
in Definition \ref{operations} but instead of using the monoidal
structures of $\FGD$ we use the corresponding structures in
$\mathbb{Z}_{2}\mbox{-}gpd.$ The negative functor
$N:\mathbb{Z}_{2}\mbox{-}gpd^{\mathbb{B}^n} \rightarrow
\mathbb{Z}_{2}\mbox{-}gpd^{\mathbb{B}^n}$ is given by
$N(F)(x,f)=N(F(x,f))$ for $F$ in
$\mathbb{Z}_{2}\mbox{-}gpd^{\mathbb{B}^n}$. The valuation map $|\
\ |: Ob(\mathbb{Z}_{2}\mbox{-}gpd^{\mathbb{B}^n}) \rightarrow
\mathbb{Q}[[x_1,...,x_n]]$ is given by $|F|= \sum_{a \in
\mathbb{N}^{n}}|F([a])|\frac{x^{a}}{a!}$.

\begin{thm}
$(\mathbb{Z}_{2}\mbox{-}gpd^{\mathbb{B}^n},N, |\ \ |)$ is a
surjective categorification of $\mathbb{Q}[[x_1,...,x_n]].$
\end{thm}

\noindent  Let us introduce examples of negative rational species.
The integral sine  $\mathrm{Si}:\mathbb{B}\rightarrow
\mathbb{Z}_{2}\mbox{-}gpd$ sends a finite set $x$ to the $\mathbb{Z}_{2}$-groupoid
$(-1)^{|x|}\mathrm{Si}(x)$ where $\mathrm{Si}(x)$ is given by
\begin{equation*}
Ob(\mathrm{Si}(x))=\left\{\begin{array}{cc} \{x\}& \ \mathrm{if}\ |x|\ \mathrm{is}\ \mathrm{odd}, \\
\emptyset & \ \mathrm{otherwise}.\ \end{array}\right.
\end{equation*}
\begin{equation*}
\mathrm{Si}(x)(x,x)=\mathbb{Z}_{|x|}\ \mathrm{for}\
|x| \ \mathrm{odd}.
\end{equation*}
\begin{prop}
\begin{equation*}
|\mathrm{Si}|=\int{\frac{\sin(x)}{x}dx}.
\end{equation*}
\end{prop}
\begin{proof}
\begin{equation*}
|\mathrm{Si}|=\sum_{n=0}^{\infty}\frac{(-1)^{n}}{2n+1}\frac{x^{2n+1}}{(2n+1)!}
=\int{\frac{\sin(x)}{x}dx}.
\end{equation*}
\end{proof}

\noindent Recall that any set $x$ may be regarded as the groupoid
whose set of objects is $x$ and with identity morphisms only.
\begin{defi}
Let $G$ be a groupoid and $n$ a positive integer. The increasing
factorial groupoid $G^{(n)}$ is given by
$$G^{(n)}=G\times(G\sqcup[1])
\times\dots\times(G\sqcup[n-1]).$$
\end{defi}
%\noindent The increasing factorial functor $( \ \
%)^{(N)}:\FGD\longrightarrow\FGD$ are given by
%$\des{\mathrm{G}\longmapsto\mathrm{G}^{(N)}}$
\noindent The set of objects of $G^{(n)}$ is
\begin{equation*}
\bigsqcup_{I\subseteq[n-1]}Ob\left(G\right) \times \left(\prod_{i\in
I}Ob\left(G\right)\right)\times\left(\prod_{i\in
[n-1]-I}[i]\right).
%=\left\{(I,f)\left|\begin{array}{c}
%\mathrm{Domain}\ \mathrm{of}\ f\ \mathrm{is}\ [n-1]\\
%f(j)\in\mathrm{G}\ \mathrm{if}\ j\in I,\ \mathrm{or}\ f(j)\in[j]\ \mathrm{if}\ j\not\in I\\
%\end{array}\right.\right\}
\end{equation*}
\noindent Thus objects in $G^{(n)}$ are
pairs $(I,f)$ such that $I\subseteq [n-1]$ and $f$ is a map with
domain $\{0,1,2,...,n-1\}$ such that $f(0)\in Ob(G),$ $f(i)\in
Ob(G)$ if $i\in I,$ and $f(i)\in[i]$ if $i\not\in I.$ Morphisms in
$G^{(n)}$ are given  by
\begin{equation*}
G^{(n)}\left((I,f);(J,g)\right)=
\left\{
\begin{array}{cc}
\emptyset & I\neq J\ \mathrm{or}\ f\mid_{I^{c}}\neq g\mid_{J^{c}},\\
{G\left(f(0),g(0) \right)\times\prod_{i\in
I}G\left(f(i),g(i)\right)} &
\mathrm{otherwise}.
\end{array}
\right.
\end{equation*}

\begin{prop}
$|{G}^{(n)}|=|G|^{(n)}$ for any finite groupoid $G$.
\end{prop}
\begin{proof}
\begin{eqnarray*}
\left|G^{(n)}\right|=
%\left|\mathrm{G}\times(\mathrm{G}\sqcup[1])
%\times\dots\times(\mathrm{G}\sqcup[n-1])\right|\\
\left|\prod_{i=0}^{n-1}(G\sqcup[i])\right|
=\prod_{i=0}^{n-1}\left|G\sqcup[i]\right|
=\prod_{i=0}^{n-1}(\left|G\right|+i)=\left|G\right|^{(n)}
\end{eqnarray*}
\end{proof}
\begin{defi}
Let $G$ be a groupoid  such that
$|G|=\frac{a}{b}\in\mathbb{Q}_{\geq 0}.$ Let
$\frac{1}{(1+X)^{\frac{a}{b}}}:\mathbb{B}\rightarrow
\mathbb{Z}_{2}\mbox{-}gpd$ be given on a finite set $x$ by
\begin{equation*}
\frac{1}{(1+X)^{\frac{a}{b}}}(x)=(-1)^{|x|}{G}^{(|x|)}.
\end{equation*}
\end{defi}
\begin{thm} For $\frac{a}{b}\in\mathbb{Q}_{\geq 0},$ we have $\left|\frac{1}
{\left(1+X\right)^{\frac{a}{b}}}\right|=\frac{1}{(1+x)^{\frac{a}{b}}}.$
\end{thm}
\begin{proof}
\begin{eqnarray*}
\left|\frac{1}{(1+X)^{\frac{a}{b}}}
\right|&=&\sum_{n=0}^{\infty}(-1)^{n}
\left(\frac{a}{b}\right)^{(n)}\frac{x^{n}}{n!}
%&=&\sum_{n=0}^{\infty}(-1)^{n}\frac{(a+b)\dots(a+b(n-1))}{b^{n}}\frac{x^{n}}{n!}\\
=\sum_{n=0}^{\infty}\left(-\frac{a}{b}\right)_{(n)}\frac{x^{n}}{n!}
=\frac{1}{(1+x)^{\frac{a}{b}}}.
\end{eqnarray*}
\end{proof}
\begin{cor}
Let $F$ be in
${{\mathbb{Z}_{2}\mbox{-}gpd_{0}}}^{\mathbb{B}^{n}}$. The species
$\frac{1}{(1+F)^{\frac{a}{b}}}:\mathbb{B}^{n}
\rightarrow
\mathbb{Z}_{2}\mbox{-}gpd$ given on a finite set $x$ by
\begin{equation*}
\frac{1}{(1+F)^{\frac{a}{b}}}(x)=\bigsqcup_{\pi\in\prod[x]}
(-1)^{|\pi|}G^{(|\pi|)}\times\prod_{b\in\pi}F(b),
\end{equation*}
is such that
\begin{equation*}
\left|\frac{1}{(1+F)^{\frac{a}{b}}}\right|=\frac{1}{(1+|F|)^{\frac{a}{b}}}.
\end{equation*}
\end{cor}
\begin{proof}
Follows from Theorem \ref{Teor6}.
\end{proof}
%\subsubsection*{The quotient groupoid}

%\noindent Recall that a set $X$ can be consider as a groupoid as
%follows $\ob(X)=X$ and morphisms in $X$ from $a$ to $b$ are given by
%$$X(a,b)=\left\{\begin{array}{cc} 1 & if\ a=b,\\
%\phi & otherwise. \end{array}\right.$$

\noindent We consider the multiplicative inverse of the valuation of a species.

\begin{thm}\label{Teor29}
\noindent Let $F$ be in ${\mathbb{Z}_{2}\mbox{-}gpd_{0}}^{\mathbb{B}^{n}}.$ The species
$(1+F)^{-1}:\mathbb{B}^{n}\rightarrow \mathbb{Z}_{2}\mbox{-}gpd$
given by
\begin{equation*}
\left({1+F}\right)^{-1}(x,f)=\bigsqcup_{\ \ \ x_{1}\sqcup\dots\sqcup
x_{n}=x}(-1)^{n}F(x_{1},f\mid_{x_{1}})\times\dots\times
F(x_{n},f\mid_{x_{n}}),
\end{equation*}
is such that
\begin{equation*}
|(1+F)^{-1}|=\frac{1}{1+|F|}.
\end{equation*}

\begin{proof}From the identities
\begin{eqnarray*}
\left({1+F}\right)^{-1}(x,f)=\bigsqcup_{\ \ \ x_{1}\sqcup\dots\sqcup
x_{n}=x}(-1)^{n}F(x_{1},f\mid_{x_{1}})\times\dots\times
F(x_{n},f\mid_{x_{n}})=\left(\sum_{n=0}^{\infty}(-1)^{n}F^{n}\right)(x,f),
\end{eqnarray*}
we conclude that
$|(1+F)^{-1}|=\sum_{n=0}^{\infty}(-1)^{n}|F|^{n}=\frac{1}{1+|F|}.$
\end{proof}
\end{thm}
%\begin{defi}
%Let $F\in\ob\left(\FGD^{\mathbb{B}^{n}}\right)$ and
%$a,b\in\mathbb{N},$ $a\neq 0.$\ The species $\frac{b}{a}F$ is given
%as follows
%$\frac{b}{a}F=\left(\bigsqcup_{i=1}^{b}\mathrm{Z}_{a}\right)F.$
%\end{defi}
%\noindent We have that
%\begin{equation*}
%\frac{1}{\left(\frac{a}{b}+F\right)}=\frac{a}{b}\frac{1}{\left(1+\frac{b}{a}F\right)}
%\end{equation*}
\begin{thm}\label{inverse}
Let $F$ be in ${\mathbb{Z}_{2}\mbox{-}gpd_{0}}^{\mathbb{B}^{n}}.$
The rational species $\frac{\mathbb{Z}_{a}^{\sqcup
b}}{1+\mathbb{Z}_{a}^{\sqcup b}F }:\mathbb{B}^{n} \rightarrow
\mathbb{Z}_{2}\mbox{-}gpd$ is such that
\begin{equation*}
\left|\frac{\mathbb{Z}_{a}^{\sqcup
b}}{1+\mathbb{Z}_{a}^{\sqcup b}F
}\right|=\frac{1}{\frac{a}{b}+|F|}.
\end{equation*}
\end{thm}
\begin{proof}
\begin{equation*}
\left|\frac{\mathbb{Z}_{a}^{\sqcup
b}}{1+\mathbb{Z}_{a}^{\sqcup b}F
}\right|=\frac{b}{a}\frac{1}{\left(1+\frac{b}{a}|F|\right)}
=\frac{1}{\frac{a}{b}+|F|}.
\end{equation*}
\end{proof}

\noindent It is worth while to pay attention to what Theorem
\ref{Teor29} and Theorem
\ref{inverse} say and  what they do not say.
Suppose that we have a species $F$ such that
$F(\emptyset)=\emptyset,$ then Theorem \ref{Teor29} constructs in
a \textit{functorial} way a species $(1+F)^{-1}$ such that
\begin{equation*}
|(1+F)^{-1}|=\frac{1}{1+|F|}.
\end{equation*}

\noindent Suppose now that we have a rational species such that
$F(\emptyset)= G$ is not the empty groupoid. Then $F=G + F_{+}$
where $F_{+}(x)=F(x)$ for $x$ a nonempty finite set  and
$F_{+}(\emptyset)=\emptyset$. Assume that $|G|=
\frac{a}{b}$, then according to Theorem
\ref{inverse} we have that

\begin{equation*}
\left|\frac{\mathbb{Z}_{a}^{\sqcup
b}}{1+\mathbb{Z}_{b}^{\sqcup a}F_{+} }\right| |F|= 1.
\end{equation*}
The species $\frac{\mathbb{Z}_{a}^{\sqcup
b}}{1+\mathbb{Z}_{a}^{\sqcup b}F_{+} }$ is constructed in a
canonical yet not functorial way from $F$. Notice the groupoid
$\mathbb{Z}_{a}^{\sqcup b}$  is such that $|\mathbb{Z}_{a}^{\sqcup
b}|=|G|^{-1},$ however the map sending $G$ into
$\mathbb{Z}_{a}^{\sqcup b}$ is not functorial.

\section{Bernoulli numbers and polynomials}
In this section we use Theorem \ref{Teor29} and Theorem
\ref{inverse} to provide a combinatorial interpretation for
Bernoulli numbers and polynomials. Let us first introduce a
generalization of the Bernoulli numbers. For $N$ a positive
integer, we let the $N$-projection map
$\pi_{N}:\mathbb{Q}[[x]]\rightarrow\mathbb{Q}[[x]]/\left(x^{N}\right)$
be given by
$\des{\pi_{N}\left(\sum_{n=0}^{\infty}f_{n}\frac{x^{n}}{n!}\right)=
\sum_{n=0}^{N-1}f_{n}\frac{x^{n}}{n!}}.$

%\begin{defi}
%Let $f\in\mathbb{Q}[[x]]$ and $N\in\mathbb{N}^{+}.$ The Taylor
%polynomial $T_{N}$ of degree $N$ is the map
%$\des{T_{N,f}(x)=\sum_{k=0}^{N}f_{k}\frac{x^{k}}{k!}}$
%\end{defi}
\begin{defi}
Let $f=\sum_{n}^{\infty}f_{n}\frac{x^n}{n!}\in\mathbb{Q}[[x]]$ and
$N\in\mathbb{N}_{+}$ be such that $f_{N}\neq 0$. The sequence
$\{B_{N,n}^{f}\}_{n=0}^{\infty}$ is called the $(f,N)$-Bernoulli
numbers sequence. It is such that
\begin{equation*}
\des{\frac{x^{N}/N!}{f(x)-\pi_{N}(f)(x)}=
\sum_{n=0}^{\infty}B_{N,n}^{f}\frac{x^{n}}{n!}}.
\end{equation*}
\end{defi}
\noindent If $N=1$ and $f(x)=e^{x}$, the sequence $B_{1,n}^{e^{x}}$
is the Bernoulli numbers sequence $B_{n}$ such that
\begin{equation*}\label{Bn,1}
\des{\frac{x}{e^{x}-1}= \sum_{n=0}^{\infty}B_{n}\frac{x^{n}}{n!}}.
\end{equation*}
\noindent If $N=2$ and $f(x)=e^{x}$ we obtain the $B_{2,n}$
Bernoulli numbers studied in  \cite{HassenNguyen}. They satisfy
\begin{equation*}\label{Bn,2}
\des{\frac{x^{2}/2!}{e^{x}-1-x}=\sum_{n=0}^{\infty}B_{2,n}\frac{x^{n}}{n!}}.
\end{equation*}
\noindent From Example \ref{Exa12} we obtain the following result.

\begin{prop}
$$|\partial{\mathbb{Z}}|=\frac{e^{x}-1}{x}.$$
\end{prop}
\noindent Theorem \rm{\ref{Teor29}} implies the following result.
\begin{thm}\label{Cor36}
The rational species $\frac{1}{\partial\mathbb{Z}}:\mathbb{B}
\rightarrow
\mathbb{Z}_{2}\mbox{-}gpd$ is such that
%$\left|\frac{1}{\partial\mathbb{Z}}\right|$
\begin{equation*}
\left|\frac{1}{\partial\mathbb{Z}}\right|=\sum_{n=0}^{\infty}B_{n}\frac{x^{n}}{n!}.
\end{equation*}
\end{thm}
%\begin{proof}
%\begin{equation*}
%\left|\frac{1}{\partial\mathbb{Z}}\right|=\frac{x}{e^{x}-1}=
%\sum_{n=0}^{\infty}B_{n}\frac{x^{n}}{n!}
%\end{equation*}
%\end{proof}
\noindent  We write $\partial{\mathbb{Z}}=1+{\partial\mathbb{Z}}_{+},$ where
${\partial\mathbb{Z}}_{+}(x)=\partial\mathbb{Z}(x)$ if
$x\neq\emptyset$ and
${\partial\mathbb{Z}}_{+}(\emptyset)=\emptyset$. We get
\begin{eqnarray*}
\frac{1}{\partial\mathbb{Z}}(x)&=&\bigsqcup_{x_{1}\sqcup\dots\sqcup
x_{k}=x}(-1)^{k}{\partial\mathbb{Z}}_{+}(x_{1})\times\dots\times
{\partial\mathbb{Z}}_{+}(x_{k})\\
&=&\bigsqcup_{x_{1}\sqcup\dots\sqcup
x_{k}=x}(-1)^{k}\mathbb{Z}(x_{1}\sqcup*_{1})\times\dots\times\mathbb{Z}(x_{k}\sqcup*_{k}).
\end{eqnarray*}
\noindent This implies that
\begin{equation*}
\frac{1}{\partial\mathbb{Z}}(x)=\bigsqcup_{x_{1}\sqcup\dots\sqcup
x_{k}=x}(-1)^{k}
\mathbb{Z}_{|x_{1}|+1}\times\dots\times\mathbb{Z}_{|x_{k}|+1}.
\end{equation*}
Next result gives a combinatorial interpretation of Bernoulli
numbers in terms of the cardinality of $\mathbb{Z}_{2}$-graded
groupoids.
\begin{cor}\label{interBerno}
\begin{equation*}
B_n=|\bigsqcup_{x_{1}\sqcup\dots\sqcup x_{k}=x}(-1)^{k}
\mathbb{Z}_{|x_{1}|+1}\times\dots\times\mathbb{Z}_{|x_{k}|+1} \ \ |.
\end{equation*}
\end{cor}
\begin{cor}
\begin{equation*}
B_{n}=\des{\sum_{a_{1}+\dots+a_{k}=n}\frac{(-1)^{k}n!}{(a_{1}+1)!\dots(a_{k}+1)!}}.
\end{equation*}
\end{cor}
\begin{proof} By Corollary \ref{interBerno}  we
have
\begin{eqnarray*}
B_{n}&=&\des{\sum_{a_{1}+\dots+a_{k}=n}
\frac{(-1)^{k}}{(a_{1}+1)\dots(a_{k}+1)}}\small{\left(\begin{array}{lll}\ \ \ \ \  n \\ \\
a_{1},\dots, a_{k}\end{array}\right)}.
%\\
%&=&\des{\sum_{a_{1}+\dots+a_{k}=n}(-1)^{k}\frac{n!}{(a_{1}+1)!\dots(a_{k}+1)!}}
\end{eqnarray*}
\end{proof}
%\txt\footnotesize{$\left(\begin{array}{lll}\ \ \ \ \  n \\ \\
%a_{1}\dots a_{k}\end{array}\right)$}
%such that
%$\partial\mathbb{Z}(x)=\mathbb{Z}(x\sqcup\{*\})$
\noindent The
decreasing factorial rational species
$\mathbb{Z}_{(N)}:\mathbb{B}\rightarrow\FGD$ is such that for each
finite set $x,$ $\ob(\mathbb{Z}_{(N)}(x))=\{x\}$ if $|x|
\geq N$ and empty otherwise. For $|x|
\geq N$ we have
$$\mathbb{Z}_{(N)}(x,x)=\mathbb{Z}_{|x|}\times\mathbb{Z}_{|x|-1}\times\dots\mathbb{Z}_{|x|-N+1}.$$
\begin{prop}
\begin{equation*}
\left|1 + N!\partial^{N}(\mathbb{Z}_{(N)})_+\right| = N!\frac{e^{x}-\pi_{N}(e^{x})}{x^{N}}.
%=\sum_{n=N}^{\infty}\frac{1}{(n)_{N}}\frac{x^{n}}{n!}
%&=&\frac{1}{N!}\int^{(N)}\left(\frac{e^{x}-\pi_{N}(e^{x})}{x^{N}/N!}\right)dx^{N}
\end{equation*}
\end{prop}
\begin{proof}
\noindent Since
\begin{equation*}
\left|\mathbb{Z}_{(N)}\right| =
\sum_{n=N}^{\infty}\frac{1}{(n)_{N}}\frac{x^{n}}{n!},
\end{equation*}
we conclude
\begin{equation*}
\des{\left| 1 + N!\partial^{N}(\mathbb{Z}_{(N)})_+\right|=
N!\sum_{n=N}^{\infty}\frac{x^{n-N}}{n!}}
=N!\frac{e^{x}-\pi_{N}(e^{x})}{x^{N}}.
\end{equation*}
\end{proof}
\noindent More generally we have the following result.
\begin{prop}
Let $F:\mathbb{B}\rightarrow \mathbb{Z}_{2}\mbox{-}gpd$ be a
rational species such that $F([N])=1$. Then
\begin{equation*}
\left|1+ N!\partial^{N}\left(F\times
\mathbb{Z}_{(N)}\right)_+\right|=N!\frac{|F|(x)-\pi_{N}(|F|)(x)}{x^{N}}.
\end{equation*}
\end{prop}
\begin{proof}
Since $$|F \times \mathbb{Z}_{(N)}|=
\sum_{n=N}^{\infty}\frac{|F([n])|}{(n)_{N}}\frac{x^{n}}{n!},$$
then $$|1 + N!\partial^{N}\left(F \times
\mathbb{Z}_{(N)}\right)_+| =
N!\sum_{n=N}^{\infty}|F([n])|\frac{x^{n-N}}{n!}=
N!\frac{|F|(x)-\pi_{N}(|F|)(x)}{x^{N}} .$$
\end{proof}
\noindent Recall that any group ${G}$ may be regarded as the groupoid
$\overline{G}$ with object $1$ and $\overline{G}(1,1)=G.$

\begin{thm}\label{bernoulligeneralizado}
Let $F:\mathbb{B}\rightarrow \mathbb{Z}_{2}\mbox{-}gpd$ be a
rational species such that $F\left([n]\right)=1$. The valuation of
the species $\left(1+ N!\partial^{N}\left(F\times
\mathbb{Z}_{(N)}\right)_+\right)^{-1}:\mathbb{B}
\rightarrow
\mathbb{Z}_{2}\mbox{-}gpd$ is the
generating series of the \ $(|F|,N)$ Bernoulli numbers. That is
\begin{equation*}
\left|\left(1+ N!\partial^{N}\left(F\times
\mathbb{Z}_{(N)}\right)_+\right)^{-1}\right|
%=\frac{x^{N}}{f(x)-\pi_{N}(f)(x)}
=\sum_{n=0}^{\infty}B_{N,n}^{|F|}\frac{x^{n}}{n!}
\end{equation*}
\end{thm}

\smallskip

\noindent Let us try to digest the meaning of Theorem
\ref{bernoulligeneralizado}. For a finite set $x$ we have
\begin{equation*}
\left(1+ N!\partial^{N}\left(F\times
\mathbb{Z}_{(N)}\right)_+\right)^{-1}(x)=
\bigsqcup_{\sqcup_{i=1}^{k}x_i=x}(-1)^{k}\prod_{i=1}^{k}N!
\partial^{N}\left(F\times\mathbb{Z}_{(N)}\right)_+(x_i).
\end{equation*}
\noindent Thus we get
\begin{equation*}
\left(1+ N!\partial^{N}\left(F\times
\mathbb{Z}_{(N)}\right)_+\right)^{-1}(x)
=
\bigsqcup_{\sqcup_{i=1}^{k}x_i=x}(-N!)^{k}
\prod_{i=1}^{k}F(x_i \sqcup [N]) \times
\prod_{i=1}^{k}\mathbb{Z}_{(N)}(x_i \sqcup [N]).
\end{equation*}

\noindent Let us now introduced a generalization of Bernoulli
polynomials.
%\section{Bernoulli polynomials}
\begin{defi}
Let $f\in\mathbb{Q}[[x]]$ and $N\in\mathbb{N}_{+}$ be such that
$f_{N}\neq 0.$ The sequence $B_{N,n}^{f}(x)$ is called the $(f,N)$
Bernoulli polynomials sequence and is such that
\begin{equation*}
\des{\sum_{n=0}^{\infty}B_{N,n}^{f}(x)\frac{y^{n}}{n!}=
\frac{\left(y^{N}/N!\right)f(xy)}{f(y)-\pi_{N}(f)(y)}}.
\end{equation*}
\end{defi}
\noindent For $N=1$ and $f(y)=e^{y},$ we obtain the Bernoulli
polynomials $B_{n}(x)$ given by
\begin{equation*}
\des{\sum_{n=0}^{\infty}B_{n}(x)\frac{y^{n}}{n!}=\frac{ye^{xy}}{e^{y}-1}}.
\end{equation*}
\noindent For $N=2$ and $f(y)=e^{y}$, we obtain the $B_{2,n}(x)$ polynomials
given by
\begin{equation*}
\des{\sum_{n=0}^{\infty}B_{2,n}(x)\frac{y^{n}}{n!}=
\frac{\left(y^{2}/2!\right)e^{xy}}{e^{y}-1-y}}.
\end{equation*}
\noindent For finite sets $a$ and $b$ we use the
notation $\mathrm{Bij}(a,b)=\{f:a\rightarrow b\mid f\ \mathrm{is}\
\mathrm{bijective}\}.$
\noindent Let $\mathrm{XZ}:\mathbb{B}^{2}\rightarrow\FGD$ be
the species such that for each pair of finite sets $(a,b)$ one has
$$\mathrm{XY}(a,b)=\left\{\begin{array}{cc} 1 & if\
|a|=|b|=1,\\ \emptyset & \mathrm{otherwise}.\end{array}
\right.$$ The valuation of $XY$ is given by $|\mathrm{XY}|=xy
\in\mathbb{Q}[[x,y]]. $

\begin{prop}\label{fin}
Let $F\in{\mathbb{Z}_{2}\mbox{-}gpd}^{\mathbb{B}^2}.$\ The
rational species $F(\mathrm{XY})$ is such that
$F(\mathrm{XY})(a,b)=F(a)\mathrm{Bij}(a,b)$  for each pair of
finite sets $a$ and $b$. Moreover $|F(XZ)|=|F|(xz)$.
\end{prop}
\begin{proof}
\begin{eqnarray*}
F(\mathrm{XY})(a,b)&=&\bigsqcup_{\pi\in\prod[x]}F(\pi)\times\prod_{b\in\pi}\mathrm{XY}(a,b)
=\bigsqcup_{\underset{\ \rm{bijective}}{f:a\rightarrow
b}}F(a)=F(a)\times\mathrm{Bij}(a,b).
\end{eqnarray*}
The identity $\left|F(\mathrm{XY})\right|=|F|(xy)$ follows from
Theorem \ref{Teor6}.
\end{proof}
\noindent Our next result provides combinatorial interpretation for
Bernoulli polynomials.

\begin{thm}
The rational species $\frac{F(\mathrm{XY})}{1+
N!\partial^{N}\left(F\times
\mathbb{Z}_{(N)}\right)_+}:\mathbb{B}
\rightarrow
\mathbb{Z}_{2}\mbox{-}gpd$ is such that
\begin{equation*}
\left|\frac{F(\mathrm{XY})}{1+ N!\partial^{N}\left(F\times
\mathbb{Z}_{(N)}\right)_+}\right|
=\sum_{n=0}^{\infty}B_{N,n}^{f}(x)\frac{y^{n}}{n!}.
\end{equation*}
\end{thm}

\noindent From Definition \ref{operations}, Theorem \ref{bernoulligeneralizado} and
Proposition \ref{fin} we see that for any pair of finite sets
$(a,b)$ the groupoid
\begin{equation*}
\frac{F(\mathrm{XY})}{1+ N!\partial^{N}\left(F\times
\mathbb{Z}_{(N)}\right)_+}(a,b)
\end{equation*}
\noindent is given by
\begin{equation*}
\bigsqcup \left( (-N!)^{k} F(c) \times
\prod_{i=1}^{k}F(x_i \sqcup [N]) \times
\prod_{i=1}^{k}\mathbb{Z}_{(N)}(x_i \sqcup [N]) \right)
\end{equation*}
\noindent where the disjoint union above runs over
subsets all $c \subseteq a$, all injective maps $i:c \rightarrow
b$, and all ordered partitions $x_1\sqcup...\sqcup x_k = b
\setminus i(c).$

\section{Euler numbers and polynomials}

In this section we provide a combinatorial interpretation for
Euler numbers and polynomials. It would be interesting  to extend
our results to the $q$-Euler numbers and polynomials discussed in
\cite{Kim}. The Euler numbers are denoted by $E_{n}$ and satisfy
\begin{equation*}
\des{\frac{2}{1+e^{x}}=\sum_{n=0}^{\infty}E_{n}\frac{x^{n}}{n!}}.
\end{equation*}
\noindent The following result should be clear.
\begin{thm}
The rational species
$\frac{1}{1+\overline{\mathbb{Z}_{2}}\mathrm{Exp}_{+}}:\mathbb{B}
\rightarrow
\mathbb{Z}_{2}\mbox{-}gpd$ is such
that
$$\left|\frac{1}{1+\overline{\mathbb{Z}_{2}}\mathrm{Exp}_{+}}\right|=\frac{2}{1+e^{x}}.$$
\end{thm}

\noindent Explicitly for any finite set $x$ we have

$$\frac{1}{1+\overline{\mathbb{Z}_{2}}\mathrm{Exp}_{+}}(x) =
\sum_{x_1 \sqcup ... \sqcup
x_k=x}(-1)^{k}(\overline{\mathbb{Z}_2})^{k}.
$$

\noindent Next we provide a combinatorial interpretation for
Euler polynomials.
\begin{defi}
The Euler polynomials sequence is given by
\begin{equation*}
\sum_{n=0}^{\infty}E_{n}(x)\frac{y^{n}}{n!}=\frac{2e^{xy}}{1+e^{y}}.
\end{equation*}
\end{defi}

\noindent We state our final result.
\begin{thm}
\begin{equation*}
\left| \frac{\mathrm{Exp}(XY)}{1+\overline{\mathbb{Z}_{2}}\mathrm{Exp}_{+}(Y)}\right|
=\sum_{n=0}^{\infty}E_{n}(x)\frac{y^{n}}{n!}.
\end{equation*}
\end{thm}
\noindent Explicitly, for any pair of finite sets $(a,b)$ we have

$$\frac{\mathrm{Exp}(XY)}{1+\overline{\mathbb{Z}_{2}}\mathrm{Exp}_{+}(Y)}(a,b)
=\bigsqcup (-1)^{k}(\overline{\mathbb{Z}_{2}})^{k},$$
\noindent where the disjoint union runs over all subsets $c \subseteq
a,$ all injective maps $i:c \rightarrow b,$ and all ordered
partitions $x_1 \sqcup ... \sqcup x_k = b\setminus i(c).$\\

\noindent It would be interesting to find combinatorial interpretation for
other known sequences of rational numbers. A step forward in that
direction has been taken in \cite{Blan}, where a combinatorial
interpretation for hypergeometric functions is provided.

\subsection* {Acknowledgment}
Our interest in the combinatorics of Bernoulli numbers arouse from
a series of lectures on \textit{Heap of Pieces} taught by Xavier
Viennot in Caracas 2002. Thanks to Mauricio Angel, Edmundo
Castillo, Allen Knutson, Jos\'{e} Mijares, Eddy Pariguan and
Domingo Quiroz. We also thank an anonymous referee for many
valuable suggestions.

\noindent hblandin@euler.ciens.ucv.ve \\
ragadiaz@gmail.com \\

\noindent Escuela de Matem\'aticas, Universidad Central de
Venezuela, Caracas 1020, Venezuela.

\end{document}